\documentclass[12pt]{article}
\usepackage{amssymb,amsmath,epsfig,amscd,eucal}

\newtheorem{thm}{Theorem}[section]
\newtheorem{defn}[thm]{Definition}
\newtheorem{lem}[thm]{Lemma}
\newtheorem{prop}[thm]{Proposition}

\newtheorem{rem}[thm]{Remark}
\newtheorem{cor}[thm]{Corollary}

\newtheorem{qu}[thm]{Question}

\newenvironment{pf}{\par\medskip\noindent{\em Proof. }}{\hfill $\square$\par\medskip}
\newenvironment{pfof}[1]{\par\medskip\noindent{\em Proof of #1. }}{\hfill $\square$\par\medskip}

\newcommand{\R}{\mathbb{R}}
\newcommand{\Z}{\mathbb{Z}}

\newcommand{\Hyp}{\mathbb{H}}

\newcommand{\lcm}{\mathrm{lcm}}

\input xy
\xyoption{all}

\title{Virtual retractions, conjugacy separability and omnipotence}
\author{Henry Wilton}
\date{$12^{\mathrm{th}}$ September 2008}

\begin{document}

\maketitle

\begin{abstract}
We use wreath products to provide criteria for a group to be conjugacy separable or omnipotent.  These criteria are in terms of virtual retractions onto cyclic subgroups.  We give two applications: a straightforward topological proof of the theorem of Stebe that infinite-order elements of Fuchsian groups (of the first type) are conjugacy distinguished, and a proof that surface groups are omnipotent.
\end{abstract}

\section{Introduction}

\begin{defn}
An element $g$ of a group $G$ is called \emph{conjugacy distinguished} if, whenever $h\in G$ is not conjugate to $g$, there exists a homomorphism $q$  to a finite group such that $q(g)$ is not conjugate to $q(h)$.  A group $G$ is called \emph{conjugacy separable} if every element of $G$ is conjugacy distinguished.
\end{defn}

The similar notion of \emph{subgroup} separability has strong connections with topology---work of Scott \cite{scott_subgroups_1978} and Stallings \cite{stallings_topology_1983} demonstrates its pleasing reformulation in terms of promoting immersions to embeddings in finite-sheeted covers.  Moreover, subgroup separability is a commensurability invariant.  In contrast, conjugacy separability is not a commensurability invariant, and does not seem to have a simple interpretation as a statement about covering spaces.  Whereas much recent work on subgroup separability relates to low-dimensional topology, the field of conjugacy separability has retained a more algebraic flavour.

In connection with the virtually Haken conjecture for hyperbolic 3-manifolds, Long and Reid made the following definitions.

\begin{defn}[\cite{long_subgroup_2008}]
A subgroup $H$ of $G$ is a \emph{virtual retract} if there exists a finite-index subgroup of $K$ such that $H$ is a retract of $K$---that is, $H\subset K$ and the inclusion map has a left inverse.
\end{defn}

\begin{defn}[\cite{long_subgroup_2008}]\label{d: LR over Z}
A group $G$ has \emph{property LR over $\Z$} if every every infinite cyclic subgroup is a virtual retract of $G$.
\end{defn}

In this paper, we provide a topological approach to proving conjugacy separability.   Lemma \ref{l: Conjugacy separability criterion} shows how a strengthening of LR over $\Z$ can be used to deduce that a group is conjugacy separable.  As an application of Lemma \ref{l: Conjugacy separability criterion}, we provide a topological proof of the theorem of Stebe \cite{stebe_conjugacy_1972} that infinite-order elements of Fuchsian groups (of the first type) are conjugacy distinguished.  It follows immediately that surface groups are conjugacy separable.  Stebe also showed in \cite{stebe_conjugacy_1972} that certain Fuchsian groups are conjugacy separable, and Fine and Rosenberger \cite{fine_conjugacy_1990} extended Stebe's result to all Fuchsian groups.

Two key steps in the proof are interesting in their own right.  Let $\alpha$ and $\beta$ be two non-homotopic curves on a surface $\Sigma$.  Proposition \ref{p: Reducing intersection number} uses Niblo's theorem that surface groups are double-coset separable \cite{niblo_separability_1992} to simplify the intersection of the elevations of $\alpha$ and $\beta$ in a finite cover.  And Proposition \ref{p: Making elevations non-homologous} produces a finite-sheeted covering in which no pair of elevations of $\alpha$ and $\beta$ are homologous.

Our motivation for studying conjugacy separability from a topological point of view is its connection with a famous open problem in geometric group theory.

\begin{qu}
Does there exist a non-residually finite hyperbolic group?
\end{qu}

Although being conjugacy separable is much stronger than being residually finite, there is a deep connection between the two.  Indeed, combinatorial Dehn filling can be used to provide a related property.

Call a set of group elements $g_1,\ldots,g_l\in G$ \emph{independent} if $i=j$ whenever $g_i$ has a conjugate that commutes with $g_j$.  Denote by $o(g)$ the order of a group element $g$.  Wise made the following definition in \cite{wise_subgroup_2000}.

\begin{defn}
A group $G$ is \emph{omnipotent} if, whenever $g_1,\ldots,g_l$ is an independent set of elements, there is an integer $K$ such that for any choice of positive integers $n_1,\ldots,n_l$, there is a homomorphism $q$ from $G$ to a finite group $Q$ such that
\[
o(q(g_i))=Kn_i
\]
for all $i$.
\end{defn}

If every hyperbolic group is residually finite then every torsion-free hyperbolic group is omnipotent, by a result of Gromov \cite{gromov_hyperbolic_1987}, Olshanskii \cite{olshanskii_residualing_1993} and Delzant \cite{delzant_sous-groupes_1996}.  Wise \cite{wise_subgroup_2000} observed that any omnipotent group that is also residually odd is conjugacy separable, and asked if there is an omnipotent group that is not conjugacy separable. Because non-residually finite hyperbolic groups seem extremely difficult to construct, one expects many hyperbolic groups to be conjugacy separable and omnipotent.

Lemma \ref{l: Omnipotence criterion} strengthens the hypotheses of Lemma \ref{l: Conjugacy separability criterion}, and enables one to deduce omnipotence. As an application, we provide a topological proof that hyperbolic surfaces are omnipotent.  This result was also obtained by Jitendra Bajpai in his Masters Thesis \cite{bajpai_omnipotence_2007}.  Indeed, we prove that such a property holds for any finite, independent set of torsion-free elements of a Fuchsian group.

A similar criterion was used by Wise \cite{wise_subgroup_2000} to prove that free groups are omnipotent.   Our techniques are slightly different to his, and when adapted to the setting of compact graphs, provide a different proof of his theorem.  To apply Lemma \ref{l: Omnipotence criterion}, we strengthen Proposition \ref{p: Making elevations non-homologous} still further: given non-homotopic curves $\alpha$ and $\beta$ on a surface, we construct a finite-sheeted cover in which any one elevation of $\alpha$ is linearly independent in homology from the set of all elevations of $\beta$.

One advantage of the techniques presented here is that they seem well adapted to dealing with finite extensions.  Although conjugacy separability, and presumably omnipotence, are not invariants of commensurability, there is no significant technical difficulty in extending our techniques from surface groups to infinite-order elements of Fuchsian groups.

Throughout this paper, if $g$ and $h$ are group elements then $g^h=h^{-1}gh$.

\subsection*{Acknowledgements}

Thanks to Alan Reid and to Martin Bridson for patient conversations.  Thanks also to Eric Katerman for some useful comments.

\section{Virtual retractions and wreath products}

For a finite-index subgroup $K$ of $G$ and an element $g\in G$, the \emph{degree} of $g$ in $K$, $\deg_K(g)$, is the minimal positive integer $n$ such that $g^n\in K$.

\begin{lem}\label{l: Conjugacy separability criterion}
Let $G$ be a group and let $a,b\in G$.  Let $K$ be a finite-index normal subgroup of $G$, let $m=\deg_K(a)$ and $n=\deg_K(b)$.
Suppose there exists a retraction $\rho:K\to\langle a^m\rangle$ with the property that
\[
\rho((b^{mn})^g)\neq a^{mn}
\]
for all $g\in G$.  Then there exists a homomorphism $\tau$ from $G$ to a finite group such that $\tau(a)$ is not conjugate to $\tau(b)$.
\end{lem}

\begin{rem}
If $g,h\in K$ then, because the image of $\rho$ is abelian, $\rho(g^h)=\rho(g)$.  Therefore, when applying Lemma \ref{l: Conjugacy separability criterion}, it suffices to check the hypotheses whenever $g$ is one of a fixed set of double-coset representatives for $K\backslash G/Z_G(b)$ (where $Z_G(b)$ is the centralizer of $b$).
\end{rem}

The idea of the proof of the lemma is to use a wreath product to promote $\rho$ to a map from $G$ to a virtually abelian group.  We therefore recall the definition of a wreath product.

Let $A$ and $B$ be groups.  Consider the group $A^B$ of set maps $B\to A$, with group operation inherited from $A$---this can be thought of as a direct sum of copies of $A$ indexed by $B$.  Then $B$ acts naturally on $A^B$ by left translation---if $\phi\in A^B$ then $\phi^b(b')=\phi(bb')$ for any $b, b'\in B$.

\begin{defn}
The \emph{wreath product} of $A$ by $B$, denoted $A\wr B$, is defined as the semidirect product
\[
A\wr B=A^B \rtimes B
\]
where the action of $B$ on $A^B$ is by left translation.  The subgroup $A^B$ is called the \emph{base} of the wreath product.
\end{defn}

Wreath products are useful because they enable us to extend homomorphisms from a normal subgroup to the whole group.  The following lemma is essentially the Krasner--Kaloujnine Theorem, which asserts that every extension is a subgroup of a wreath product.

\begin{lem}
Suppose $K\lhd G$ with $Q\cong G/K$ and let $f:K\to H$ be a homomorphism.  Then there is a homomorphism
\[
\hat{f}:G\to H\wr Q
\]
with the property that, whenever $k\in K$, $\hat{f}(k)$ is in the base of the wreath product and, for any $q\in Q$, there is $g_q\in G$ such that
\[
\hat{f}(k)(q)=f(k^{g_q}).
\]
\end{lem}
\begin{pf}
Denote by $\eta$ the quotient map $G\to Q$, and for each $q\in Q$ fix a corresponding coset representative $g_q\in G$, so that $\eta(g_q)=q$.  For any $g\in G$ and $q\in Q$ set
\[
\phi_g(q)=g_q^{-1}gg_{\eta(g)^{-1}q}.
\]
Note that $\eta(\phi_g(q))=1$, so $\phi_g\in K^Q$ and therefore the composition $f\circ\phi_g$ makes sense as an element of $H^Q$.   Now define $\hat{f}:G\to H\wr Q$ by
\[
\hat{f}(g)=(f\circ\phi_g)\eta(g).
\]
Let us check that this is a homomorphism.  For any $g,g'\in G$,
\[
\hat{f}(g)\hat{f}(g')=(f\circ\phi_g)\eta(g)(f\circ\phi_{g'})\eta(g')=(f\circ\phi_g)(f\circ\phi_{g'})^{\eta(g)^{-1}}\eta(gg').
\]
For any $q\in Q$,
\begin{eqnarray*}
(f\circ\phi_g)(f\circ\phi_{g'})^{\eta(g)^{-1}}(q)&=&(f\circ\phi_g)(q)(f\circ\phi_{g'})(\eta(g)^{-1}q)\\
&=&f(g_q^{-1}gg_{\eta(g)^{-1}q})f(g_{\eta(g)^{-1}q}^{-1}g'g_{\eta(g')^{-1}\eta(g)^{-1}q})\\
&=&f(g_q^{-1}gg'g_{\eta(gg')^{-1}q})\\
&=&f\circ\phi_{gg'}(q)
\end{eqnarray*}
and so
\[
\hat{f}(g)\hat{f}(g')=(f\circ\phi_{gg'})\eta(gg')=\hat{f}(gg')
\]
as required.

That $\hat{f}(k)(q)=f(k^{g_q})$ for $k\in K$ and $q\in Q$ is immediate from the construction.
\end{pf}

We will need to test when elements of wreath products are conjugate.

\begin{rem}\label{r: Conjugacy in wreath products}
Suppose $A$ is abelian.  If $\phi_1, \phi_2\in A^B$ are conjugate in $A\wr B$ then for some $b\in B$ we have $\phi_1^b=\phi_2$.  As $B$ acts on $A^B$ by permuting the factors, it follows that the sets
\[
\{\phi_1(b)\mid b\in B\}
\]
and
\[
\{\phi_2(b)\mid b\in B\}
\]
are equal.
\end{rem}

We now have the tools to prove Lemma \ref{l: Conjugacy separability criterion}.

\begin{pfof}{Lemma \ref{l: Conjugacy separability criterion}}
Let $Q=K\backslash G$. Consider the composition
\[
\sigma:K\stackrel{\rho}{\to} \Z \to \Z/N
\]
for $N$ large enough that $\sigma((b^{mn})^g)\neq \sigma(a^{mn})$ for all $g\in G$.  Let $\tau$ be the extension of $\sigma$ to a map $G\to (\Z/N)\wr Q$.

Suppose $\tau(a)$ and $\tau(b)$ are conjugate.  Then $\tau(a^{mn})$ and $\tau(b^{mn})$ are conjugate.  But $a^{mn}$ and $b^{mn}$ are both in $K$, so it follows from Remark \ref{r: Conjugacy in wreath products} that for some $q\in Q$,
\[
\sigma(a^{mn})=\sigma((b^{mn})^{g_q}),
\]
a contradiction.
\end{pfof}

To provide a criterion for omnipotence, we need to analyse the order of $\tau(b)$.  Let $\sigma':K\to (\Z/N)^Q$ be the restriction of $\tau$ to $K$.  Thinking of $\rho$ as a map $K\to\Z$, let $d_b=\gcd\{\rho((b^n)^g)|g\in G\}$.  Then, as long as $N$ is sufficiently large, $\sigma'(b^n)$ is a primitive element of $(\Z/N)^Q$ multiplied by $d_b$.  So $o(\sigma'(b^n))=\lcm(N,d_b)/d_b$  (adopting the convention that $0/0=1$).   As $o(q(b))=n$, it follows that
\[
o(\tau(b))=n\frac{\lcm(N,d_b)}{d_b}.
\]
As $\rho(a^m)=1$, we have that $d_a=1$ and therefore $o(\tau(a))=mN$.  Therefore
\[
\frac{o(\tau(a))}{o(\tau(b))}=\frac{mNd_b}{n\lcm(N,d_b)}
\]
and so, unless $d_b=0$,
\[
\frac{m}{n}\leq\frac{o(\tau(a))}{o(\tau(b))}\leq\frac{md_b}{n}.
\]
To prove omnipotence, we shall need the ratio of the orders of $\tau(a)$ and $\tau(b)$ to vary unrestrictedly with $N$; we therefore need $d_b=0$.

Bearing this in mind, Lemma \ref{l: Omnipotence criterion} provides a criterion to prove omnipotence for groups.

\begin{lem}\label{l: Omnipotence criterion}
Let $G$ be a group and let $\{a_1,\ldots,a_l\}$ be an independent set of elements.  Suppose that there is a finite-index normal subgroup $K\lhd G$, with $m_i=\deg_{K}(a_i)$, and suppose further that for each $i$ there exists a retraction $\rho_i:K\to\langle a_i^{m_i}\rangle$ with the property that whenever $j\neq i$,
\[
\rho_i((a_j^{m_j})^g)= 1
\]
for all $g\in G$.  Then for any choice of positive integers $p_1,\ldots,p_l$ there exists a homomorphism $\eta$ from $G$ to a finite group such that
\[
o(\eta(a_i))=p_i\prod_j m_j
\]
for all $i$.
\end{lem}
\begin{pf}
Applying the wreath product construction above to each $\rho_i$, we see that for any choices of positive integer $N_i$, for each $i$ there is a homomorphism $\sigma_i$ from $G$ to a finite group $Q_i$ such that
\[
o(\sigma_i(a_i))=N_im_i
\]
whereas
\[
o(\sigma_i(a_j))=m_j
\]
whenever $j\neq i$.  For each $i$ let
\[
N_i=p_i\prod_{j\neq i}m_{j}
\]
and let $\sigma_i:G\to Q_i$ be the resulting homomorphism to a finite group.  Now
\[
\eta=\prod_i\sigma_i:G\to \prod_i Q_i
\]
is easily seen to be as required.
\end{pf}

\section{Conjugacy separability}

We shall use the ideas of section 1 to give a topological proof of the theorem of Stebe \cite{stebe_conjugacy_1972} that infinite-order elements of Fuchsian groups of the first type are conjugacy distinguished.

\begin{defn}
A \emph{Fuchsian group of the first type} is a discrete subgroup of $PSL_2(\R)$ of finite covolume.
\end{defn}

In what follows, for brevity's sake we shall simply refer to Fuchsian groups, when we really mean Fuchsian groups of the first type.

We can think of a Fuchsian group $\Gamma$ as acting properly discontinuously on the hyperbolic plane $\Hyp^2$, such that the quotient is a cone-type 2-orbifold.  By Selberg's Lemma $\Gamma$ has a torsion-free normal subgroup $\Gamma_0$ of finite index, and the quotient of $\Hyp^2$ by $\Gamma_0$ is a surface $\Sigma$.  Because $\Gamma$ and hence $\Gamma_0$ is finitely generated we can restrict out attention to a compact subsurface (which we also denote $\Sigma$), possibly with boundary, whose fundamental group is $\Gamma_0$.   The proof that Fuchsian groups are conjugacy separable (and, later, omnipotent) proceeds by analysing closed curves on the compact surface $\Sigma$.

Given a closed curve $\gamma$ on $\Sigma$, the homology class (with $\Z$ coefficients) of $\gamma$ is denoted $[\gamma]$.   A closed curve $\gamma$ that is not null-homotopic is called \emph{primitive} if it is simple and $[\gamma]$ is primitive in $H_1(\Sigma)$.  (So either $\gamma$ is non-separating or $\gamma$ is boundary parallel and $\Sigma$ has more than one boundary component.)  If $\alpha,\beta$ are closed curves on $\Sigma$ then $i(\alpha,\beta)$ is the geometric intersection number of $\alpha$ and $\beta$.

Throughout the following we shall use the language of \emph{elevations}, which for the purposes of this paper we define as follows.

\begin{defn}
Let $\gamma$ be a closed curve on $\Sigma$ and let $\Sigma'\to\Sigma$ be a covering map.  If $g\in\pi_1(\Sigma)$ is freely homotopic to $\gamma$ and $n=\deg_{\pi_1(\Sigma')}(g)$ then, by standard covering space theory, $g^n$ lifts to a curve $\gamma'$ on $\Sigma'$.  Such a curve $\gamma'$, defined up to free homotopy on $\Sigma'$, is called an \emph{elevation} of $\gamma$ to $\Sigma'$.

As free homotopy corresponds to conjugation in the fundamental group, one can equivalently think of $\gamma'$ as the $\pi_1(\Sigma')$-conjugacy class of $g^n$.\footnote{Wise uses a slightly different definition of elevations extensively in \cite{wise_subgroup_2000}.  His definition is better adapted to more general contexts, although for our purposes the two definitions coincide.  His definition of elevations is naturally in bijection with the set of double cosets $\pi_1(\Sigma')\backslash\pi_1(\Sigma)/\langle\gamma\rangle$, whereas the definition we use here is naturally in bijection with the set of double cosets $\pi_1(\Sigma')\backslash\pi_1(\Sigma)/Z_{\pi_1(\Sigma)}(\gamma)$}
\end{defn}

Scott famously proved that surface groups are subgroup separable \cite{scott_subgroups_1978}.   We shall use Niblo's extension of this result.  He proved that surface groups are double-coset separable.

\begin{thm}[\cite{niblo_separability_1992}]
Suppose that $H, H'$ are finitely generated subgroups of $\pi_1(\Sigma)$ and $\gamma\notin HH'$.  There exists a finite-index subgroup $K\subset \pi_1(\Gamma)$ such that $H\subset K$ but $\gamma\notin KH'$.
\end{thm}

Scott's Theorem has the well known consequence that any closed curve on a surface can be lifted to a simple curve in a finite-sheeted cover.  We shall use Niblo's Theorem to simplify the intersections of a pair of curves.  It is also well known that any separating simple closed curve can be lifted to be non-separating in a finite cover.  It will be important later to understand this fact, so we recall the details here.

\begin{lem}\label{l: Making separating curves non-separating}
Let $\gamma$ be a separating curve on a surface $\Sigma$, let $\Sigma'\to\Sigma$ be a finite-sheeted abelian covering and let $\gamma'$ be the lift of $\gamma$ to $\Sigma'$.  If $\gamma'$ is non-separating then the restriction of the covering map $\Sigma'\to\Sigma$ to some component of $\Sigma'\smallsetminus\gamma'$ is a homeomorphism.
\end{lem}
\begin{pf}
Let $A$ be the abelian group of covering transformations of $\Sigma'$ and let $\Sigma'_1$ and $\Sigma'_2$ be the components of $\Sigma'\smallsetminus\gamma'$.  For a contradiction, suppose that there are non-trivial $a_1,a_2\in A$ such that $a_i\gamma'\subset\Sigma'_i$ for $i=1,2$.  Without loss of generality, we can assume that both $a_i\gamma'$ are closest to $\gamma'$, in the sense that there is a path from $\gamma'$ to $a_i\gamma'$ that does not cross any other translates of $\gamma'$.

An easy induction shows that $a^n_1\gamma'\subset\Sigma'_1$ for any positive integer $n$ and hence, because $A$ is finite, $a_1^{-1}\gamma'\subset\Sigma'_1$.  It follows that $a_1\Sigma'_2\subset\Sigma'_1$.  Likewise, $a_2\Sigma'_1\subset\Sigma'_2$.  Therefore
\[
a_1a_2\gamma'\subset \Sigma'_1
\]
and, symmetrically,
\[
a_2a_1\gamma'\subset\Sigma'_2.
\]
Because $A$ is abelian, $a_1a_2=a_2a_1$ so $a_1a_2$ maps $\gamma'$ to itself.  As $A$ acts freely on the elevations of $\gamma$ it follows that $a_1a_2=1$ and $a_2=a_1^{-1}$.  But we have already seen that $a_1^{-1}\gamma'\subset\Sigma'_1$, so it follows that $a_1$ and $a_2$ fix $\gamma'$, a contradiction.

We conclude that, without loss of generality, every translate of $\gamma'$ by $A$ is contained in $\Sigma'_2$.  Hence the restriction of the covering map to $\Sigma'_1$ is a homeomorphism.
\end{pf}

We are now ready to apply Niblo's Theorem to simplify the intersections of a pair of closed curves.

\begin{prop}\label{p: Reducing intersection number}
Let $\Sigma$ be a compact surface and let $\alpha, \beta$ be closed curves in $\Sigma$ that are not null-homotopic.  There is a finite-sheeted covering $\Sigma'\to\Sigma$ and a primitive elevation $\alpha'$ of $\alpha$ to $\Sigma'$ such that, for any elevation $\beta'$ of $\beta$ to $\Sigma'$:
\begin{enumerate}
\item $\beta'$ is primitive; and
\item $i(\alpha',\beta')\leq 1$.
\end{enumerate}
\end{prop}
\begin{pf}
By passing to a double cover if necessary, it is easy to ensure that if $\Sigma$ has non-empty boundary then it has at least two boundary components.   It follows from Scott's Theorem that there is an orientable finite-sheeted covering $\hat{\Sigma}$ of $\Sigma$ to which $\alpha$ has a lift $\hat{\alpha}$ that is a simple closed curve.

If $\alpha$ is boundary parallel then, as $\hat{\Sigma}$ has at least two boundary components, $\hat{\alpha}$ is primitive.  If $\alpha$ is non-separating then it is easy to apply Lemma \ref{l: Making separating curves non-separating} and pass to a double cover so that $\hat{\alpha}$ is non-separating and hence primitive.  Likewise we can pass to a further finite-sheeted cover to ensure that some elevation of $\beta$ is also primitive. Replacing $\hat{\Sigma}$ by a normal covering, we can ensure that every elevation of $\alpha$ and $\beta$ is primitive.

We now need to simplify the intersections of an elevation $\hat{\alpha}$ of $\alpha$ with the elevations of $\beta$.  If $\hat{\alpha}$ is boundary parallel then any curve in $\hat{\Sigma}$ can be homotoped off $\hat{\alpha}$.  We shall therefore concentrate on the case in which $\hat{\alpha}$ is non-separating.

Let $\{\hat{\beta}_j\}$ be the set of all elevations of $\beta$ to $\hat{\Sigma}$.  The proof is by induction on the quantity
\[
c(\hat{\alpha},\hat{\Sigma})=\sum_j \max(i(\hat{\alpha},\hat{\beta}_j)-1,0).
\]
Clearly, if $c(\hat{\alpha},\hat{\Sigma})=0$ then $i(\hat{\alpha},\hat{\beta_j})\leq 1$ for all $j$.

If $c(\hat{\alpha},\hat{\Sigma})>0$ then without loss of generality $i(\hat{\alpha},\hat{\beta}_1)>1$.  After modifying $\hat{\alpha}$ and $\hat{\beta}_1$ by a homotopy we may assume that $|\hat{\alpha}\cap\hat{\beta}_1|=i(\hat{\alpha},\hat{\beta}_1)$.  Fix a basepoint $\hat{x}\in\hat{\alpha}\cap\hat{\beta_1}$.  Then $\hat{\beta}_1$ is homotopic (respecting the basepoint) to a concatenation $\gamma\delta$ where $i(\hat{\alpha},\gamma)\leq 1$ and $i(\hat{\alpha},\delta)<i(\hat{\alpha},\hat{\beta}_1)$.

Whenever $g\in\langle\hat{\alpha}\rangle\langle\hat{\beta}_1\rangle$, either $g\in\langle\hat{\alpha}\rangle$ or $g=\hat{\alpha}^m\hat{\beta}_1^n$ for some integer $m$ and some $n\neq 0$, so
\[
i(\hat{\alpha},g)=i(\hat{\alpha},\hat{\beta}_1^n)>1.
\]
But $\gamma\notin\langle\hat{\alpha}\rangle$, as otherwise a homotopy would reduce $|\hat{\alpha}\cap\hat{\beta}_1|$, so $\gamma\notin \langle\hat{\alpha}\rangle\langle\hat{\beta}_1\rangle$.  Therefore by Niblo's Theorem there exists a (based) finite-sheeted covering
\[
(\bar{\Sigma},\bar{x})\to(\hat{\Sigma},\hat{x})
\]
such that $\hat{\alpha}$ lifts to a (based) loop $\bar{\alpha}$ on $\bar{\Sigma}$ but $\gamma\notin\pi_1(\bar{\Sigma})\langle\hat{\beta}_1\rangle$.

We aim to show that $c(\bar{\alpha},\bar{\Sigma})<c(\hat{\alpha},\hat{\Sigma})$.  Fix $j$, and consider the set $\{\bar{\beta}_{j,k}\}$ of elevations of $\hat{\beta}_j$ to $\bar{\Sigma}$.  Then it is clear that
\[
\sum_k i(\bar{\alpha},\bar{\beta}_{j,k})= i(\hat{\alpha},\hat{\beta}_j).
\]
Therefore $c(\bar{\alpha},\bar{\Sigma})\leq c(\hat{\alpha},\hat{\Sigma})$.

Let $\bar{\beta}_1$ be the elevation of $\hat{\beta}_1$ to $\bar{\Sigma}$ that covers the based loop $\hat{\beta}_1$ and let $\bar{\beta}'_1$ be the elevation that covers $\hat{\beta}_1^\gamma$.  Both $\bar{\beta}_1$ and $\bar{\beta}'_1$ intersect $\bar{\alpha}$ non-trivially.  But $\bar{\beta}_1$ and $\bar{\beta}'_1$ are not freely homotopic in $\pi_1(\bar{\Sigma})$.  This shows that $\hat{\beta}_1$ has at least two elevations to $\bar{\Sigma}$ that intersect $\bar{\alpha}$ non-trivially.  It follows that
\[
c(\bar{\alpha},\bar{\Sigma})<c(\hat{\alpha},\hat{\Sigma})
\]
as required.

By induction, there exists a finite-sheeted covering $\Sigma'\to\Sigma$ with an elevation $\alpha'$ of $\alpha$, such that $c(\alpha',\Sigma')=0$.   This is the required covering.
\end{pf}

To apply the results of section 1, we need to separate elevations of $\alpha$ and $\beta$ in an abelian quotient, and hence in homology.  Our first step is the following observation.

\begin{rem}\label{r: Making disjoint curves non-homologous}
If $\alpha$ and $\beta$ are disjoint, primitive curves in $\Sigma$ that are homologous but not homotopic then there exists a 2-sheeted cover $\Sigma'$ of $\Sigma$ to which $\alpha$ and $\beta$ both lift, such that no pair of elevations $\alpha'$ and $\beta'$ of $\alpha$ and $\beta$ respectively are homologous in $\Sigma'$.  Let us describe this covering.

As usual, the easier case is when $\alpha$ is boundary parallel.  In this case $\beta$ is also boundary parallel, and $\Sigma$ only has these two boundary components.  A double cover $\Sigma'$ with four boundary components is as required.

In the case when $\alpha$ and $\beta$ are non-separating, note that $\Sigma\smallsetminus\alpha\cup\beta$ has two components $\Sigma_1$ and $\Sigma_2$, each of which is a surface of positive genus with two boundary components.  Each $\Sigma_i$ has a double cover $\Sigma'_i$ with four boundary components.  Gluing $\Sigma'_1$ and $\Sigma'_2$ suitably creates $\Sigma'$.
\end{rem}

This remark enables us to improve Proposition \ref{p: Reducing intersection number} to ensure that the elevations of $\alpha$ and $\beta$ differ in homology.

\begin{prop}\label{p: Making elevations non-homologous}
Let $\Sigma$ be a surface and let $\alpha, \beta$ be closed curves in $\Sigma$ that are not freely homotopic.  There is a finite-sheeted normal covering $\tilde{\Sigma}\to\Sigma$ such that, for any elevations $\tilde{\alpha}$ and $\tilde{\beta}$ of $\alpha$ and $\beta$ respectively to $\tilde{\Sigma}$:
\begin{enumerate}
\item $\tilde{\alpha}$ and $\tilde{\beta}$ are primitive; and
\item $[\tilde{\alpha}]\neq \pm[\tilde{\beta}]$ in $H_1(\tilde{\Sigma})$.
\end{enumerate}
\end{prop}
\begin{pf}
Let $\Sigma'$ be the covering space of $\Sigma$ given by Proposition \ref{p: Reducing intersection number}, with the provided elevation $\alpha'$ of $\alpha$.   Let $\{\beta'_j\}$ be the set of all elevations of $\beta$ to $\Sigma'$.  For each $j$, define a covering space $\Sigma'_j$ as follows.  If $\alpha'$ is not homologous to $\beta'_j$ or $(\beta'_j)^{-1}$ in $\Sigma'$ then $\Sigma'_j=\Sigma'$.  Otherwise, $\alpha'$ and $(\beta'_j)^{\pm 1}$ are (homotopic to) disjoint primitive curves in $\Sigma'$ that are homologous but not homotopic, so by Remark \ref{r: Making disjoint curves non-homologous} we can take $\Sigma'_j$ to be the double cover of $\Sigma'$ in which no pair of elevations of $\alpha'$ and $(\beta'_j)^{\pm 1}$ are homologous.  Define $\bar{\Sigma}$ by
\[
\pi_1(\bar{\Sigma})=\bigcap_j\pi_1(\Sigma'_j).
\]
Then whenever $\bar{\alpha}$ is a lift of $\alpha'$ to $\bar{\Sigma}$ and $\bar{\beta}$ is an elevation of $\beta$ to $\bar{\Sigma}$,  $\bar{\beta}$ is also an elevation of some intermediate $\beta'_j$ and hence is homologous to neither $\alpha'$ nor its inverse.

Finally, let $\tilde{\Sigma}$ be the covering space whose fundamental group is the intersection of all conjugates of $\pi_1(\bar{\Sigma})$ in $\pi_1(\Sigma)$.  This is the covering required.  For, after composing the covering map $\tilde{\Sigma}\to\bar{\Sigma}$ with a deck transformation, any pair $\tilde{\alpha}$ and $\tilde{\beta}$ of elevations (of $\alpha$ and $\beta$ respectively) cover a lift $\bar{\alpha}$ of $\alpha'$ and some elevation $\bar{\beta}$ of $\beta$.  Hence $\tilde{\alpha}$ and $\tilde{\beta}^{\pm 1}$ are not homologous.
\end{pf}

Combining this proposition with the results of section 1, we are now in a position to prove that Fuchsian groups are conjugacy separable.  Like Stebe, we shall start by proving that elements of infinite order are conjugacy distinguished.

\begin{thm}[\cite{stebe_conjugacy_1972}]\label{t: Infinite-order elements are conjugacy distinguished}
If $\Gamma$ is a Fuchsian group and $a\in\Gamma$ is of infinite order then $a$ is conjugacy distinguished.
\end{thm}
\begin{pf}
By Selberg's Lemma, $\Gamma$ has a torsion-free normal subgroup of finite index, which can be taken to be $\pi_1(\Sigma)$ for some compact surface $\Sigma$, possibly with boundary.  Suppose that $b\in\Gamma$ is not conjugate to $a$.  Let $p=\deg_{\pi_1(\Sigma)}(a)$ and $q=\deg_{\pi_1(\Sigma)}(b)$, and represent $a^p$ by a closed curve $\alpha$ on $\Sigma$.  Fix representatives $g_1,\ldots,g_n$ for the set of double cosets $\pi_1(\Sigma)\backslash \Gamma/Z_\Gamma(b)$ and, for each $i$, let $\beta_i$ be a closed curve on $\Sigma$ that represents $(b^q)^{g_i}$.

If $b$ is of finite order then every $\beta_i$ is null-homotopic.  By Proposition \ref{p: Reducing intersection number} there is a finite-sheeted cover $\Sigma'$ to which $\alpha$ has a primitive elevation $\alpha'$.  Because $\alpha'$ is primitive in $H_1(\Sigma')$ there is a retraction $\rho:\pi_1(\Sigma')\to\langle\alpha'\rangle$, and
\[
\rho(\beta'_i)=1
\]
whenever $\beta'_i$ is an elevation of some $\beta_i$ to $\Sigma'$.  It now follows from Lemma \ref{l: Conjugacy separability criterion} that there is a homomorphism to a finite group under which the images of $a$ and $b$ are non-conjugate.

We can therefore assume that $b$ is of infinite order.  For each $i$, let $\tilde{\Sigma}_i$ be the finite-sheeted covering provided by Proposition \ref{p: Making elevations non-homologous}, in which no elevation of $\alpha$ is homologous to an elevation of $\beta_i^{\pm 1}$.  Now let $\tilde{\Sigma}$ be the covering defined by
\[
\pi_1(\tilde{\Sigma})=\bigcap_i\pi_1(\tilde{\Sigma}_i).
\]
Intersecting $\pi_1(\tilde{\Sigma})$ with its conjugates, we can assume furthermore that $\pi_1(\tilde{\Sigma})$ is a normal subgroup of $\Gamma$.  Let $m=\deg_{\pi_1(\tilde{\Sigma})}(a)$ and let $\tilde{\alpha}$ be a closed curve representing $a^m$ in $\pi_1(\tilde{\Sigma})$.

If $n=\deg_{\pi_1(\tilde{\Sigma})}(b)$ then any $\Gamma$-conjugate of $b^n$ is conjugate in $\pi_1(\tilde{\Sigma})$ to some elevation of some $\beta_i$.  Hence, to satisfy the hypotheses of Lemma \ref{l: Conjugacy separability criterion}, we must construct a retraction $\rho:\pi_1(\tilde{\Sigma})\to\langle\tilde{\alpha}\rangle$ such that $\rho(\tilde{\beta}_i^m)\neq\tilde{\alpha}^n$ whenever $\tilde{\beta}_i$ is an elevation of some $\beta_i$.

Suppose that $n[\tilde{\beta}_i]=m[\tilde{\alpha}]$ in $H_1(\tilde{\Sigma})$.  Then, as $\tilde{\alpha}$ and $\tilde{\beta}_i$ are both primitive in homology, $m=n$ and $[\tilde{\alpha}]=\pm[\tilde{\beta}_i]$.  But this contradicts the properties of $\tilde{\Sigma}$.  Therefore, choosing a suitable projection $H_1(\tilde{\Sigma})\to\langle\tilde{\alpha}\rangle$ and defining $\rho$ to be the concatenation
\[
\rho:\pi_1(\tilde{\Sigma})\to H_1(\tilde{\Sigma})\to \langle \tilde{\alpha}\rangle
\]
we obtain a retraction that satisfies the hypotheses of Lemma \ref{l: Conjugacy separability criterion}.  It follows that there exists a homomorphism $\sigma$ from $\pi_1(\Sigma)$ to a finite group such that $\sigma(a)$ and $\sigma(b)$ are not conjugate, as required.
\end{pf}

Because surface groups are torsion-free, it follows immediately that hyperbolic surface groups are conjugacy separable.  Indeed, the proof did not use the hyperbolic structure on $\Sigma$, and applies just as well to Euclidean surfaces.

\begin{cor}[\cite{fine_conjugacy_1990}]\label{c: Fuchsian groups are conjugacy separable}
Surface groups are conjugacy separable.
\end{cor}

\begin{rem}\label{r: Free groups are conjugacy separable}
That free groups are conjugacy separable is a special case of Theorem \ref{t: Infinite-order elements are conjugacy distinguished}.  Indeed, one can provide a very quick proof of this fact using Lemma \ref{l: Conjugacy separability criterion} and the topology of graphs.  The idea is to pass to a finite-sheeted cover in which every elevation of a pair of non-conjugate elements is a simple closed curve.  It is then clear that they differ in homology.
\end{rem}

\section{Omnipotence}

In this section we improve further upon the results of the previous section to prove omnipotence.  To apply Lemma \ref{l: Omnipotence criterion}, we need a criterion to ensure that elevations are linearly independent in homology.  The next proposition provides this.

\begin{prop}\label{p: Making elevations independent in homology}
Let $\Sigma$ be a hyperbolic surface and let $\alpha,\beta_1,\ldots,\beta_l$ be a collection of primitive curves with the property that no $\beta_i$ is homologous to $\alpha$ or its inverse.  Then there exists a finite-sheeted cyclic covering space $\check{\Sigma}$ such that any elevation $\check{\alpha}$ of $\alpha$ is linearly independent in homology of the set of all elevations of all the $\beta_i$.  Hence, there exists a retraction
\[
\rho:\pi_1(\check{\Sigma})\to\langle\check{\alpha}\rangle
\]
such that $\rho(\check{\beta}_i)=1$ whenever $\check{\beta}_i$ is an elevation of some $\beta_i$.
\end{prop}
\begin{pf}
Consider first the case when $\alpha$ is non-separating.  Let $\gamma_1$ be a non-separating simple closed curve with $i(\alpha,\gamma_1)=1$.  The commutator $\delta=[\alpha,\gamma_1]$ is a separating simple closed curve, because $\Sigma$ is hyperbolic. Let $\gamma_2$ be a primitive curve in the component of $\Sigma\smallsetminus\delta$ that does not contain $\gamma_1$.  Neither $\gamma_1$ nor $\gamma_2$ is homologous to $\alpha$ or its inverse.

By the hypotheses on homology there exists a homomorphism $\pi_1(\Sigma)\to\Z$ that kills $\alpha$ but kills none of the $\beta_i$ and neither of the $\gamma_k$.  Consider the composition
\[
\psi:\pi_1(\Sigma)\to\Z\to\Z/p
\]
where $p$ is a prime large enough that all the $\psi(\beta_i)$ and both the $\psi(\gamma_k)$ are non-zero.  Let $\iota:\check{\Sigma}\to\Sigma$ be the covering map with $\pi_1(\check{\Sigma})=\ker\psi$.  Then $\alpha$ has $p$ elevations to $\check{\Sigma}$, which we denote by $\check{\alpha}_1,\ldots,\check{\alpha}_p$, whereas each $\beta_i$ only has one elevation, denoted $\check{\beta}_i$.

Without loss of generality we can take each $\check{\alpha}_j$ to be a lift of $\alpha^{\gamma_1^{1-j}}$.  In particular, $\check{\alpha}_1\check{\alpha}^{-1}_2=\delta$ in $\pi_1(\Sigma)$.  But because neither $\gamma_1$ nor $\gamma_2$ lift to $\check{\Sigma}$, it follows from Lemma \ref{l: Making separating curves non-separating} that $\delta$ lifts to a non-separating curve $\check{\delta}$ in $\check{\Sigma}$.  Hence
\[
[\check{\alpha}_1]-[\check{\alpha}_2]=[\check{\delta}]\neq 0
\]
so the elevations $\check{\alpha}_1$ and $\check{\alpha}_2$ are not homologous.

Associated to the covering map $\iota$ there is the induced map $\iota_*:H_1(\check{\Sigma})\to H_1(\Sigma)$ and the transfer map $\tau^*:H_1(\Sigma)\to H_1(\check{\Sigma})$, which maps a curve to the sum of its elevations.  The composition $\iota_*\circ\tau^*$ is equal to multiplication by $p$.  We shall use the transfer map to show that $[\check{\alpha}_1]$ is linearly independent from the $[\check{\beta}_i]$.

Suppose that there are constants $\kappa$ and $\lambda_i$ such that the equation
\[
\kappa[\check{\alpha}_1]=\sum_i\lambda_i[\check{\beta}_i]
\]
holds in $H_1(\check{\Sigma})$.  Applying $\iota_*$ to both sides gives
\[
\kappa[\alpha]=p\sum_i\lambda_i[\beta_i]
\]
so $\kappa=p\kappa'$ for some $\kappa'$.  Dividing by $p$ and applying $\tau^*$ to both sides this becomes
\[
\kappa'\sum_j[\check{\alpha}_j]=\sum_i\lambda_i[\check{\beta}_i]=\kappa[\check{\alpha}_1].
\]
The action of $\Z/p$ on $\check{\Sigma}$ is transitive on the set of elevations of $\alpha$, so unless $\kappa=0$ it follows that
\[
[\check{\alpha}_1]=[\check{\alpha}_j]
\]
for any $j$.  But this contradicts our previous observation that $\check{\alpha}_1$ and $\check{\alpha}_2$ are not homologous, so $\kappa=0$.  This completes the proof when $\alpha$ is non-separating.

Suppose now that $\alpha$ is boundary-parallel.  The proof proceeds similarly to the previous case.  As before, let $\psi:\pi_1(\Sigma)\to\Z/p$ be a homomorphism that kills $\alpha$ but kills none of the $\beta_i$, and let $\check{\Sigma}$ be the covering corresponding to $\ker\psi$.  Now, because $p\geq 2$, we see that $\check{\Sigma}$ has at least four boundary components, and hence no pair of elevations of $\alpha$ is homologous.  The same argument using the transfer map again implies that an elevation $\check{\alpha}_1$ of $\alpha$ is linearly independent from the set of elevations of the $\beta_i$.
\end{pf}

Let $\alpha_1,\ldots,\alpha_l$ be an independent collection of elements of $\pi_1(\Sigma)$, viewed as closed curves on $\Sigma$. We can apply Proposition \ref{p: Making elevations independent in homology} pairwise to the $\alpha_i$.

\begin{rem}\label{r: Making elevations of independent sets independent}
Fix $i\neq j$.  By Proposition \ref{p: Making elevations non-homologous} there exists a finite-sheeted cover $\tilde{\Sigma}_{i,j}$ of $\Sigma$ in which, whenever $\tilde{\alpha}_i$ and $\tilde{\alpha}_j$ are elevations of $\alpha_i$ and $\alpha_j$ respectively:
\begin{enumerate}
\item $\tilde{\alpha}_i$ is a non-separating simple closed curve; and
\item $\tilde{\alpha}_i$ and $\tilde{\alpha}_j^{\pm 1}$ are not homologous.
\end{enumerate}
Let $\tilde{\Sigma}$ be the finite-sheeted cover of $\Sigma$ such that
\[
\pi_1(\tilde{\Sigma})=\bigcap_{i\neq j} \pi_1(\tilde{\Sigma}_{i,j}).
\]
Then $\tilde{\Sigma}$ has the property that every elevation of an $\tilde{\alpha}_i$ is a non-separating simple closed curve, and no pair of elevations of $\tilde{\alpha}_i$ and $\tilde{\alpha}_j^{\pm 1}$ are homologous when $i\neq j$.

Fix $i$ and an elevation $\tilde{\alpha}_i$ of $\alpha_i$.  Let $\{\tilde{\alpha}_{j,k}|i\neq j\}$ be the set of all elevations of all $\alpha_j$ (where $i\neq j$) to $\tilde{\Sigma}$.  Applying Proposition \ref{p: Making elevations independent in homology} to $\tilde{\alpha}_i$ together with the set $\{\tilde{\alpha}_{j,k}|i\neq j\}$, it follows that there exists a finite-sheeted normal covering $\check{\Sigma}_i\to\tilde{\Sigma}$ with the property that any elevation $\check{\alpha}_i$ of $\tilde{\alpha}_i$ is linearly independent from the set of all elevations $\{\check{\alpha}_{j,k}|i
\neq j\}$ of all the $\tilde{\alpha}_{j,k}$ (with $i\neq j$).

Let $\check{\Sigma}$ be the cover of $\Sigma$ such that
\[
\pi_1(\Sigma)=\bigcap_i\pi_1(\check{\Sigma}_i).
\]
This has the required property for all $i$---that is, whenever $\check{\alpha}_i$ is an elevation of $\alpha_i$ to $\check{\Sigma}$, there is a retraction  $\rho_i:\check{\Sigma}\to\langle\check{\alpha}_i\rangle$ such that $\rho_i(\check{\alpha}_j)=1$ for any elevation $\check{\alpha}_j$ of any $\alpha_j$ with $i\neq j$. Furthermore, by intersecting $\pi_1(\check{\Sigma})$ with all its conjugates, we can assume that $\check{\Sigma}\to\Sigma$ is a regular covering.
\end{rem}

We can now prove our main theorem.

\begin{thm}\label{t: Omnipotence for Fuchsian groups}
Let $\Gamma$ be a Fuchsian group (of the first kind) and let $a_1,\ldots,a_l$ be an independent set of elements of infinite order in $\Gamma$.  Then there is a constant $K$ such that, whenever $p_1,\ldots,p_l$ are positive integers, there is a homomorphism $\eta$ from $\Gamma$ to a finite group such that
\[
o(\eta(a_i))=Kp_i
\]
for each $i$.
\end{thm}
\begin{pf}
By Selberg's Lemma, $\Gamma$ has a torsion-free normal subgroup of finite index, which we can think of as $\pi_1(\Sigma)$ where $\Sigma$ is a  hyperbolic surface.  For each $i$, let $m_i=\deg_{\pi_1(\Sigma)}(a_i)$ and fix a set of representatives $\{g_{i,k}\}$ for $\pi_1(\Sigma)\backslash \Gamma/Z_\Gamma(a_i)$.  For each $i$ and $k$, let $\alpha_{i,k}$ be a closed curve on $\Sigma$ representing $(a_i^{m_i})^{g_{i,k}}$.

The set of $\alpha_{i,k}$ correspond to an independent set of elements of $\pi_1(\Sigma)$, so by Remark \ref{r: Making elevations of independent sets independent} there is a finite-sheeted covering $\check{\Sigma}\to\Sigma$ such that any elevation $\check{\alpha}_{i,k}$ of $\alpha_{i,k}$ is linearly independent in homology from the set of homology classes of all elevations of the set $\{\alpha_{j,k'}|(i,k)\neq(j,k')\}$.  It follows that there exists a retraction $\rho_i:\pi_1(\check{\Sigma})\to\langle\check{\alpha}_{i,k}\rangle$ such that
\[
\rho_i(\check{\alpha}_{j,k'})=1
\]
whenever $j\neq i$ (indeed, whenever $(j,k')\neq (i,k)$).  It follows that $\rho_i$ satisfies the hypotheses of Lemma \ref{l: Omnipotence criterion} and hence, setting $K=\prod_i m_i$, there exists a map $\eta$ from $\Gamma$ to a finite group such that
\[
o(\eta(a_i))=Kp_i
\]
for each $i$, as required.
\end{pf}

It follows immediately that surface groups are omnipotent.

\begin{cor}[Bajpai \cite{bajpai_omnipotence_2007}]\label{c: Omnipotence for surface groups}
If $\Sigma$ is a compact hyperbolic surface (possibly with boundary) then $\pi_1(\Sigma)$ is omnipotent.
\end{cor}

The theorem of Wise that free groups are omnipotent \cite{wise_subgroup_2000} is a corollary.  His proof of this fact uses a similar criterion to that of Lemma \ref{l: Omnipotence criterion}, although he does not use wreath products to achieve it.  He then proceeds to construct the required retraction explicitly, using a refinement of Stallings' proof of Marshall Hall's Theorem \cite{stallings_topology_1983}.  The methods of this paper provide a simple alternative proof of omnipotence for free groups.  The ideas of Proposition \ref{p: Making elevations independent in homology} work just as well in the context of graphs, and omnipotence follows immediately from Lemma \ref{l: Omnipotence criterion}.

\bibliographystyle{plain}

\bigskip\bigskip\centerline{\textbf{Author's address}}

\smallskip\begin{center}\begin{tabular}{l}%
Department of Mathematics\\
1 University Station C1200\\
Austin, TX 78712-0257\\
USA\\
{\texttt{henry.wilton@math.utexas.edu}}
\end{tabular}\end{center}

\end{document}